\documentclass[11pt]{article}
\usepackage{latexsym,euscript,amsmath,amssymb,amsbsy,amsfonts,amsthm,amsopn,amstext,amsxtra,amscd}
\usepackage{epsfig}
\usepackage{graphics}
\usepackage{url}

\usepackage[english]{babel}

\theoremstyle{plain}
\newtheorem{theorem}{Theorem}

\newtheorem{corollary}{Corollary}
\newtheorem{proposition}{Proposition}

\theoremstyle{definition}
\newtheorem{example}{Example}
\newtheorem{definition}{Definition}

\newtheorem{remark}{Remark}

\newcommand{\N}{{\mathbb N}}
\newcommand{\R}{{\mathbb R}}

\newcommand{\Z}{{\mathbb Z}}

\newcommand{\D}{\text{$\mathcal{D}$}}

\begin{document}
\centerline{\bf  MEASURABLE SEQUENCES}

\vskip1cm

\centerline{Milan Pa\v st\'eka and Robert Tichy}

\vskip1cm

\centerline{\bf ABSTRACT}

{\it The paper deals with the distribution functions of sequences
with respect to asymptotic density and measure density. Furthermore also polyadicly continuous sequences and
their extension to random variables are studied.}

\section{Introduction}

In the first part we study sequences having an asymptotic
distribution function in the sense of Schoenberg [Sch].
The connection between independence and statistical independence is
established in case of continuous distribution functions.

 Later we develop relations between distribution
functions of sequences and distribution functions of random
variables. We study statistical independence and independence in the
sense of probability theory. In the last part we transfer some
probabilistic limit laws to certain types of deterministic
sequences. This makes heavily use of methods developed in [CQ].

The ''general'' notion of uniform distribution was introduced by Hermann Weyl (1916) in
his famous paper [WEY]: a sequence $\{v(n)\}, v(n) \in
[0,1)$ is {\it uniformly distributed} if and only if for every $x
\in [0,1)$
$$
\lim_{N \to \infty} \frac{1}{N} |\{n \le N; v(n) < x\}| =x,
$$
where $|A|$ denotes he cardinality of the set $A$. This can be
equivalently formulated using the notion of asymptotic density.
Let $\N$ be the set of positive integers. We say that a set $A
\subset \N$ has an {\it asymptotic density} if and only if the limit
$$
\lim_{N \to \infty}\frac{|A \cap [1, N)|}{N}:= d(A)
$$
exists, and in this case the value $d(A)$ is called the {\it
asymptotic density} of $A$. Let $\D$ denote the system of all
subsets of $\N$ having an asymptotic density. Then a sequence
$\{v(n)\}, v(n) \in [0,1)$ is uniformly distributed if and only if
for every $x \in [0,1)$ the set $\{n \le N; v(n) < x\}$ belongs to
$\D$ and $d(\{n \le N; v(n) < x\}) =x$. Schoenberg
[Sch] generalized this notion as follows: we say that a sequence
$\{v(n)\}, v(n) \in [0,1)$ has an {\it asymptotic distribution
function} if and only if for each real number $x$ the set $\{n \le
N; v(n) < x\}$ belongs to $\D$. In this case the function $F(x) =
d(\{n \le N; v(n) < x\})$ is called the {\it asymptotic distribution
function} of the sequence $\{v(n)\}$.

Our aim is to study distribution functions of sequences. The
following statement is useful in this context.

\begin{proposition}
\label{APP} If $F$ is a non decreasing function defined on the real line
then for each real numbers $x_1, x_2$ - the points of continuity of
$F$-we have that for every $\varepsilon>0$ there exist two
continuous function $g, g_1$ such that
$$
g \le \mathcal{X}_{[x_1,x_2]} \le g_1
$$
and
$$
\int_{-\infty}^\infty (g_1(x) - g(x)) < \varepsilon,
$$

$\mathcal{X}_{[x_1,x_2]}$ denoting the indicator function of the interval $[x_{1}, x_{2}]$.
\end{proposition}
The proof follows from a standard procedure, see [KN] page 54.\\

Another important notion of uniform distribution was introduce by Niven [NIV].
A sequence of positive
integers $k=\{k_n\}$ is called {\it uniformly distributed in} $\Z$ if
and only for each $m \in \N, r \in \Z$ we have that $\{n \in \N; k_n
\equiv r \mod{m}\} \in \D$ and $d(\{n \in \N; k_n \equiv r
\mod{m}\}) = \frac{1}{m}$. In sections 9 and 10 we will use this concept to prove structural
properties concerning measurable sequences.

\section{\bf Mean value, dispersion and Buck measurability}

Let $v =\{v(n)\}$ be a sequence of real numbers. Set
$$
E_N(v)= \frac{1}{N}\sum_{n=1}^N v(n)
$$
for $N =1,2,3,\dots$ and
$$
\underline{E}(v)=\liminf_{N \to \infty}E_N(v),
\overline{E}(v)=\limsup_{N \to \infty}E_N(v).
$$
\begin{definition}
\label{meanvalue} If $\underline{E}(v) =\overline{E}(v):=E(v)$ we
say that $v$ has a {\it mean value} and the number $E(v)$ will be
called the {\it mean value of $v$}.
\end{definition}

Clearly we have
\begin{proposition}
\label{lin} If sequences $v, w$ have mean values then for all
numbers $a, b$ the sequence $av+bw$ has a mean value and
$$
E(av+bw)=aE(v)+bE(w).
$$
\end{proposition}

\begin{proposition} If $v$ is bounded sequence with elements in the
interval $[a,b]$ and having an asymptotic distribution function $F,$ then $v$
has a mean value and
$$
E(v)=\int_a^bxdF(x).
$$
\end{proposition}

\begin{definition}
\label{dispersion} We say that a sequence $v$ has a dispersion if
$v$ has a mean value and the sequence $(v-E(v))^2$ has a mean value;
in this case the number
$$
D^2(v)=E((v-E(v))^2)
$$
is called the {\it dispersion} of $v$.
\end{definition}

If a bounded sequence $v$ has a dispersion then
$$
D^2(v)=\lim_{N \to \infty} \frac{1}{N}\sum_{n=1}^N (v(n)-E(v))^2.
$$

In the following we introduce Buck measurability, weak measurability and weak distribution functions.

R. C. Buck [BUC] constructed a measure
density via covering of sets by arithmetic progressions. Denote
$$
r+(m)= \{r+jm; j=0,1,2,\dots \}
$$
for $r =0,1,2, \dots$ and $m \in \N$. Then $r+(m)$ belongs to $\D$
and $d(r+(m)) =\frac{1}{m}$. If $S \subset \N$ then the value
$$
\mu^\ast(S) = \inf \Big\{ \sum_{j=1}^k \frac{1}{m_k}; S \subset
\bigcup_{j=1}^k r_j+(m_j) \Big\}
$$
is called {\it Buck's measure density} of the set $S$.

The sets from the system
$$
\D_\mu =\{S \subset \N; \mu^\ast(S)+\mu^\ast(\N \setminus S)=1 \}
$$
are called {\it Buck measurable}.

The following trivial fact will be useful for us, (see [PAS3], page
39) :
\begin{proposition}
\label{propI}  a) $\D_\mu$ is an algebra of sets, the restriction
$\mu=\mu^\ast|_{\D_\mu}$ is a finitely additive probability measure
on $\D_\mu$. \newline b) $\D_\mu \subset \D$ and $d(S)=\mu(S)$ for
every $S \in \D_\mu$. \newline c) A set $S \subset \N$ belongs to
$\D_\mu$ if and only if for each $\varepsilon > 0$ sets $S_1,
S_2 \in \D_\mu$ exist such that $S_1 \subset S \subset S_2$ and $\mu(S_2)
- \mu(S_1) < \varepsilon$.
\end{proposition}

We say that a sequence of real numbers $\{v(n) \}$ is {\it Buck
measurable} if and only if for every real number $x$ the set $\{n
\in \N; v(n) < x\}$ belongs to $\D_\mu$. In this case the function
$$
F(x) = \mu(\{n \in \N; v(n) < x\})
$$
is called {\it Buck's distribution function} (for short B-d.f.) of $\{v(n)\}$.

A Buck measurable sequence is called {\it Buck uniformly distributed (for short B- u.d.)
} if and only if its Buck distribution function $F(x)$ satisfies
\begin{equation}
\label{bud} F(x)=0,\text{for}\; x<0, F(x)=x,\text{for}\; x \in [0,1], F(x)=1,\text{for}\; x>1.
\end{equation}

 Proposition \ref{propI} implies
\begin{proposition} { Each Buck measurable sequence of real numbers has an
asymptotic distribution function which coincides with its Buck
distribution function.}
\end{proposition}

\begin{definition}
\label{weakmesur} A real valued sequence $\{v(n)\}$ is called {\it
weakly Buck
 measurable} if and only if the sets $\{n \in \N; v(n)<x\}$ are Buck
 measurable excluding at most a countable set of real numbers
 $x$. In this case the function
 $$
 F(x)=\{n \in \N; v(n)<x\}
 $$
 defined on the real line excluding at most a countable set is called a  {\it weak Buck distribution function} of
$\{v(n)\}$.
\end{definition}

The standard procedure yields a variant of Chebyshev's inequality:

\begin{proposition}
\label{chebyshev} If a bounded sequence $v$ has weak distribution
function then for each $\varepsilon>0$ we have
$$
\overline{d}(\{n \in \N; |v(n)-E(v)|
> \varepsilon\}) \le \frac{D^2(v)}{\varepsilon^2}.
$$
\end{proposition}

Using these concepts we obtain the following

\begin{proposition} If $v$ is a bounded sequence having a weak
distribution function and $D^2(v)=0$ then there exists a set $A \in
\D$ such that $d(A)=1$ and $\lim_A v(n)=E(v),$ were the limit is taken a along the set $A.$
\end{proposition}

\begin{definition} Let $v, w$ be sequences having weak asymptotic
distribution functions. Suppose moreover that the sequence $vw$ has a
mean value. The value
$$
\rho(v,w)= \frac{|E(vw)-E(v)E(w)|}{D(v)D(w)}
$$
will be called the {\it correlation coefficient } of the sequences $v,
w$.
\end{definition}

\begin{definition}
\label{defcor}We say that the sequences $v, w$ are {\it correlated}
if and only if such values $\alpha,  \beta$ exist that $\lim_A w(n)
- \alpha v(n) - \beta = 0$ for some set $A$ from $\D$ such that
$d(A) =1,$ where the limit is taken along the set A.
\end{definition}

In [P-T] the following result is proved

\begin{proposition} The sequences v, w are correlated if and only if
$vw$ has a mean value and $\rho(v,w) = 1$. In this case for
$\alpha, \beta$ from Definition \ref{defcor} we have
$$
\alpha= \frac{E(vw)-E(v)E(w)}{D^2(v)}, \beta = E(w)-\alpha E(v).
$$
\end{proposition}

This has the following implication:

\begin{corollary} If $v, w$ are sequences uniformly distributed
modulo 1, then they are correlated if and only if $E(vw)=\frac{1}{3}$
or $E(vw)=\frac{1}{6}$.
\end{corollary}

 \section{Independent sequences}

 In the book [Ra] the following notion is defined :

\begin{definition} Two bounded real valued sequences $v, w$ are
called {\it statistically independent} if and only if for every
functions $g, g_1$, that are continuous on a closed interval containing the
elements of both sequences
$$
\lim_{N\to \infty} E_N(g(v))E_N(g_1(w))- E_N(g(v)g_1(w))=0.
$$
\end{definition}

Properties of statistical independent sequences are studied in
various papers. For a survey we refer to the monograph [SP].

\begin{definition} Two sets $S, S_1 \in \D$ are called {\it
independent} if and only if $S\cap S_1 \in \D$ and $d(S\cap S_1) =
d(S) d(S_1)$. Bounded sequences $v, w$ having asymptotic
distribution functions are called {\it independent} if and only if
for arbitrary intervals $I, I_1$ the sets $\{n \in \N; v(n)\in I_{1}\}$
and $\{n \in \N; w(n) \in I_1\}$ are independent.
\end{definition}

We shall prove

\begin{theorem}
\label{instat} Let $v, w$ be bounded sequences having continuous
asymptotic distribution functions. Then these sequences are
independent if and only if they are statistically independent.
\end{theorem}

We start with the following

\begin{proposition}
\label{UNKO} Let $v_k, w_k, (k \in \N) $ be two systems of sequences of
elements from a certain closed interval $[a,b]$. Suppose that for each
$k \in \N$ the sequences $v_k, w_k$ are statistically independent.
If $v_k$ converges uniformly to $v$ and $w_k$ converges uniformly to
$w$ then the sequences $v,w$ are statistically independent.
\end{proposition}

{\bf Proof.} Let $g,g_1$ be continuous functions defined on a closed
interval containing the elements of both sequences. Then these
functions are uniformly continuous. Thus $g(v_k)$ converges
uniformly to $g(v)$, $g_1(w_k)$ converges uniformly to $g(w)$ and
$g(v_k)g_1(w_k)$ converges uniformly to $g(v)g_1(w)$. Hence for given
$\varepsilon>0$ there exists $k$ with
$$
|E_N(g(v_k))-E_N(g(v))|<\varepsilon,
|E_N(g(w_k))-E_N(g(w))|<\varepsilon
$$
and
$$
|E_N(g(v_k)g_1(w_k))-E_N(g(v)g_1(w))| < \varepsilon.
$$
Moreover there exists $N_0$ such that for $N \ge N_0$ we have
$$
|E_N(g(v_k))E_N(w_k)-E_N(g(v_k)g_1(w_k))| < \varepsilon.
$$
From the first inequalities we derive
$$
|E_N(g(v_k))E_N(w_k)- E_N(g(v))E_N(g_1(w))|< 2M\varepsilon,
$$
where $M$ is an upper bound of $|v|, |w|$. This yields for $N \ge
N_0$
$$
|E_N(g(v)g_1(w))- E_N(g(v))E_N(g_1(w))|< 2M\varepsilon+\varepsilon.
$$
\qed

\begin{definition}
\label{simpl} If $S_1, \dots, S_k$ are disjoint sets belonging to
$\D$ and $c_{1},\ldots c_{k} \in \mathbb{R}$ then the sequence $s$ defined by
$$
s(n)= \sum_{j=1}^k c_j\mathcal{X}_{S_j}(n), n \in \N
$$
is called a {\it simple
sequence}.
\end{definition}

It is easy to check :
\begin{proposition} If $s$ is a simple sequence then $s$ has a mean value and
$$
E(s)= \sum_{j=1}^k c_jd(S_j)
.$$
\end{proposition}

This leads to the following consequence:

\begin{proposition} Let $s=\sum_{j=1}^k c_j\mathcal{X}_{S_j}, r=\sum_{j=1}^\ell
r_j\mathcal{X}_{R_j}$ be such simple sequences that the sets $S_j,
R_m$, are independent for $j=1,\dots,k, m=1, \dots, \ell.$ Then they
are statistically independent.
\end{proposition}

\begin{proposition} If $v, w$ are bounded independent sequences then
they are statistically independent.
\end{proposition}

{\bf Proof.} Let the values of $v, w$ be contained in the interval $[a,b]$.
Consider for $k \in \N$ the partition of $[a,b]$ into disjoint
subintervals $I_j, j=1,\dots, m$ such that $|I_j|<\frac{1}{k},
j=1,\dots,m$. Then the sets
$$
S_j=\{n\in \N; v(n)\in I_j\}, R_i=\{n\in \N; w(n)\in I_i\}, 1\le i,
j \le m
$$
are independent. Thus the simple sequences
$$s_k=\sum_{j=1}^m c_j
\mathcal{X}_{S_j}, r_k=\sum_{j=1}^m c_j \mathcal{X}_{R_j}, c_j \in
I_j, j=1,\dots, m
$$ are statistically independent. Since $|s_k(n)-v(n)|\le \frac{1}{k}$
and $|r_k(n)-w(n)|\le \frac{1}{k}$ for $n \in \N$ we obtain that $s_k$
converges uniformly to $v$ and $r_k$ converges uniformly to $w$. Thus
 due to Proposition \ref{UNKO} $v$ and $w$ are
statistically independent. \qed

{\bf Proof of the second implication of  Theorem \ref{instat}}.
Consider statistically independent bounded sequences $v, w$ having
continuous asymptotic distribution functions $F, F_1$ respectively.
Let $I_1 =[x_1,x_2], I_2=[y_1,y_2]$. Since  $F, F_1,$ are continuous,
Proposition \ref{APP} guarantees that for $\varepsilon>0$ there
exist positive continuous functions $f,f_1, g, g_1$ satisfying
\begin{equation}
\label{I1} f \le \mathcal{X}_{I_1} \le f_1, g \le \mathcal{X}_{I_2}
\le g_1
\end{equation}
and
\begin{equation} \label{I2} \int_a^b (f_1(x)-f(x)) dF(x) <
\varepsilon, \int_a^b (g_1(x)-g(x)) dF_1(x) < \varepsilon.
\end{equation}
From (\ref{I1}) we derive
$$
E_N(f(v)g(w))\le E_N(\mathcal{X}_{I_1}(v)\mathcal{X}_{I_2}(w)) \le
E_N(f_1(v)g_1(w)),
$$
moreover
$$
E_N(f(v))E_N(g(w))\le
E_N(\mathcal{X}_{I_1}(v))E_N(\mathcal{X}_{I_2}(w)) \le
E_N(f_1(v))E_N(g_1(w)).
$$
If $N \to \infty$ we obtain for
$\overline{E}=\overline{E}(\mathcal{X}_{I_1}(v)\mathcal{X}_{I_2}(w))$
and
$\underline{E}=\underline{E}(\mathcal{X}_{I_1}(v)\mathcal{X}_{I_2}(w))$
the inequalities
$$
\int^a_bf(x)dF(x)\int^a_bg(x)dF_1(x) \le \overline{E}\le
\int^a_bf_1(x)dF(x)\int^a_bg_1(x)dF_1(x)
$$
and
$$
\int^a_bf(x)dF(x)\int^a_bg(x)dF_1(x) \le \underline{E}\le
\int^a_bf_1(x)dF(x)\int^a_bg_1(x)dF_1(x)
.$$ Set $S_1=\{n \in \N; v(n) \in I_1\}, S_2=\{n \in \N; w(n)
\in I_2\}$. Then
$$
\lim_{N\to \infty}
E_N(\mathcal{X}_{I_1}(v))E_N(\mathcal{X}_{I_2}(w)) =d(S_1)d(S_2).
$$
This yields
$$
|\overline{E}-d(S_1)d(S_2)| \le H\varepsilon ,
|\underline{E}-d(S_1)d(S_2)| \le H\varepsilon,
$$
where $H$ is suitable constant. Since $\varepsilon>$ is arbitrary we
get $\overline{E}= \underline{E}=d(S_1)d(S_2)$. If we consider that
$E=\overline{E}=\underline{E}=d(S_1\cap S_2)$ the assertion follows.
\qed

An immediate consequence of the definition is the following:
\begin{proposition}
\label{vect} Let $v, w$ be bounded sequences having continuous
asymptotic distributions $F, F_1,$ respectively. Suppose that these
sequences are independent. Then for any intervals $I_1=[x_1, x_2],
I_2=[y_1, y_2]$ the set $S= \{n \in \N; (v(n),w(n))\in I_1 \times
I_2\}$ belongs to $\D$ and
$$
d(S)=(F(x_2)-F(x_1))(F_1(y_1)-F_1(y_2)).
$$
\end{proposition}

Using the above notation the standard method yields:

\begin{proposition} Let $A$ be a Riemann Stjeltjes measurable set with
respect product measure $F\times F_1$ then the set $R=\{n\in \N;
(v(n), w(n)) \in A\}$ belongs to $\D$ and
$$
d(R) = \int \int_R dF(t_1)dF_1(t_2).
$$
\end{proposition}

Furthermore the following theorem holds (with the above notation).

\begin{theorem}
\label{plus}
The sequence $v+w$ has an asymptotic distribution function $F_2$ given by
$$
F_2(x)=\int\int_{\{(t_1, t_2); t_1 + t_2 \le x\}} dF(t_1)dF(t_1).
$$
\end{theorem}

This leads after some calculation to:

\begin{corollary} If $v, w$ are two
independent uniformly distributed sequences then the sequence $v+w$
has the distribution function $G$ where $G(x)=0, x\le 0,
G(x)=\frac{x^2}{2}, x \in [0,1], G(x)=2x -\frac{x^2}{2} - 1, x \in
[1,2], G(x)=1, x>2$.
\end{corollary}

In following the more general notion of independence will be useful:

\begin{definition} If $v_1, \dots, v_k$ are bounded sequences having
asymptotic distribution functions then they are called {\it
independent} if and only if for all intervals $I_1, \dots, I_k$ the set $S=\{n\in \N; v_j(n)\in I_j, j=1,\dots, k\}$
belongs to $\D$ and
$$
d(S)=\prod_{j=1}^k d(\{n \in \N; v_j(n) \in I_j\}).
$$
These sequences are called {\it statistically independent} if and
only if
$$
\lim_{N \to \infty}E_N(g_1(v_1)\dots g_k(v_k))-E_N(g_1(v_1))\dots
E_N(g_k(v_k)) =0
$$
for any functions $g_1, \dots, g_k$ continuous on closed intervals
containing all elements of the given sequences.
\end{definition}

\begin{theorem}
\label{GENIND} Let $v_1, \dots, v_k$ be bounded sequences having
continuous asymptotic distribution functions. Then they are
independent if and only if they are statistically independent.
\end{theorem}

\begin{proposition} Let $v_1, v_2, v_3$ be bounded sequences having asymptotic distribution
functions. If these sequences are independent then $v_1+v_2, v_3$ are
independent, too.
\end{proposition}

From this we derive as above:

\begin{theorem} If $v_1, \dots, v_k$ are independent bounded
sequences with continuous distribution functions, having the same
mean value $E$ and the same dispersion $D^2.$ Then
$$
d\Big(\Big\{n \in \N; \Big|\frac{v_1+\dots+v_k}{k}-E\Big| \ge
\varepsilon \Big\}\Big) \le \frac{D^2}{n\varepsilon^2}.
$$
\end{theorem}

\section{Polyadicly continuous sequences}

Denote by $\Omega$ the compact metric ring of polyadic integers, (see
[N], [N1], [PAS5], which is the
completion of $\N$ with respect to the polyadic metric
\begin{equation}
\label{metric} \mathfrak{d}(a,b) =\sum_{n=1}^\infty
\frac{\psi_n(a-b)}{2^n},
\end{equation}
where $\psi_n(x)=0$ if $n$ devides $x$ and $\psi_n(x)=1$ otherwise. For sequences $\{v(n)\}$ we shall
use two synonymous expressions: sequences or arithmetic
functions. A sequence $\{v(n)\}$ is called {\it
polyadicly continuous} (for short: {\it$p$-continuous)} if and only if for each $\varepsilon >0$ there is $m
\in \N$ such that
$$
\forall a, b \in \N; a\equiv b \pmod{m} \Rightarrow |v(a)-v(b)| <
\varepsilon.
$$

In [PAS2] it is proved:
\begin{proposition}
\label{wk} { Let $\{v(n)\}$ be a $p$-continuous sequence of
elements of $[0,1]$. Suppose that $F$ is a continuous function
defined on $[0,1]$. Then $\{v(n)\}$ is Buck measurable with B-d.f. $F$ if and only if
$$
\lim_{N \to \infty} \frac{1}{N}\sum_{n=1}^N h(v(n)) = \int_0^1
h(x)dF(x)
$$
for each continuous real valued function $h$ defined on $[0,1]$.}
\end{proposition}

This implies the following
\begin{proposition}
\label{erd} {If a $p$-continuous sequence of elements in
$[0,1]$ has a continuous asymptotic distribution function then it is Buck
measurable and its B-d.f. coincides with its
asymptotic distribution function.}
\end{proposition}

The next result is due to P. Erd\H{o}s [Er], see also [PAS3], p.32.
In the following we use the notation
$$A(v(n),I) = \mid \{n\in \mathbb{N}; v(n)\in I\}\mid,$$ for sequences $v(n)\in I \subset [0,1].$

\begin{theorem}

{ Suppose that $f$ is a non-negative additive arithmetic
  function such that for every prime $p$ we have $f(p)=f(p^k), k=1,2,3,
  \dots$,
  and for distinct primes \index{prime} $p_1,p_2$ we have
$f(p_1) \neq f(p_2)$. Assume that the infinite series $\sum_{p
\index{prime}}\frac{f(p)}{p}$ (running over the primes) converges. Then for every interval
$I$, there holds $A(\{f(n)\},I) \in \D$. Moreover, in this case, the
function $$g(x) = d(A(\{f(n)\}, [-\infty, x))$$ is continuous on the real
line.}
\end{theorem}

\begin{corollary} { Let $f$ be a non-negative additive arithmetic
function  such that for every prime $p$ we have $f(p)=f(p^k),
k=1,2,3,  \dots$, for different primes $f(p_1) \neq f(p_2)$ and the series $\sum_p f(p)$ converges.
Then the sequence $\{f(n)\}$ is Buck measurable with continuous
Buck distribution function.}
\end{corollary}

{\bf Proof.} Let $N \in \N$. If $n_1 \equiv n_2 \pmod{N!}$ then
$n_1,n_2$ contain the same primes smaller than $N$ in canonical
decomposition and so in this case
$$
|f(n_1)-f(n_2)| \le 2\sum_{p > N} f(p)
.$$

Thus the convergence of $\sum_p f(p)$ provides that $\{f(n)\}$ is a
$p$-continuous sequence. This condition yields also the
convergence of $\sum_p \frac{f(p)}{p},$ and the assertion follows.\qed\\

It is easy to check that each $p$-continuous sequence of
real numbers is uniformly continuous with respect to the polyadic metric
$\mathfrak{d},$ and so each  $p$-continuous sequence of real
numbers $\{v(n)\}$ can be extended in the natural way to a real
valued continuous function $\tilde{v}$ defined on $\Omega$ such that
$$
\tilde{v}(\alpha) = \lim_{j\to \infty}v(n_j),
$$
where $\{n_j\}$ is a sequence of positive integers such that $n_j \to
\alpha$ for $j \to \infty$ with respect the polyadic metric. The
compact ring $\Omega$ is equipped with Haar probability measure $P$
and so the function $\tilde{v}$ can be considered as random variable
on the probability space $(\Omega, P)$. As usually $h$ is a
random variable on $\Omega$ and we denote $E(h)=\int h dP$, the mean
value of $h$.

Let $m \in \N$ and $s = 0,1, \dots m-1$. Put
$$
s+m\Omega = \{s+m\alpha; \alpha \in \Omega\}.
$$
The ring $\Omega$ can be represented as disjoint union
$$
\Omega = \bigcup_{s=0}^{m-1} s +m\Omega.
$$
(see [N], [N1]).  Thus for the Haar probability measure $P$ we
have
\begin{equation}
\label{mesure} P(s +m\Omega)= \frac{1}{m}
\end{equation}
for $m \in \N, s=0, \dots m-1$.

In [PAS4] it is proven that
\begin{equation}
\label{clo}
 \mu^\ast(S) = P(cl(S))
\end{equation}
for each $S \subset \N$, where $cl(S)$ denote the topological
closure of $S$ in $\Omega$.

\begin{example}
\label{vdc} Let $\{Q_k\}$ be an increasing sequence of integers such
that $Q_0=1$ and $Q_k \;\text{devides}\; Q_{k+1}, k=1,2,3,\dots$. Each positive
integer $n$ can be uniquely represented in the form
$$
n = a_0+ a_1Q_1+ ...+ a_kQ_k
,$$
where $a_j < \frac{Q_{j+1}}{Q_j}, j=1,\dots ,k$. To this $n$ we
associate an element $\gamma(n)$ in the unit interval of the form
$$
\gamma(n)= \frac{a_0}{Q_1}+\dots + \frac{a_k}{Q_{k+1}}.
$$
The sequence $\{\gamma(n)\}$ is known as van der Corput sequence
in base $\{Q_k\}$ and in [PAS2] it is proved that it is Buck uniformly distributed and $p$-continuous.
\end{example}

The following characterization allows us to apply results of
probability theory to the distribution of $p$-continuous
sequences:

\begin{theorem}
\label{thm1}
Let $\{v(n)\}$ be a $p$-continuous sequence and
$F$ a continuous real valued function defined on the real line. Then the
following the statements are equivalent: \end{theorem}
  (i) $F$ {\it is the
distribution function of the random variable} $\tilde{v}$.

 (ii)
$\{v(n)\}$ {\it is a Buck measurable sequence and $F$ is its  B-d.f.}.

 (iii) {\it For each real number $x$ we have}
$$
\mu^\ast(\{n \in \N; v(n) < x \}) = F(x).
$$

{\bf Proof.}(i) $\Rightarrow$ (ii). The continuity of $F$ yields
\begin{equation}
\label{e2} P(\tilde{v} < x) = F(x)=P(\tilde{v} \le x)
\end{equation}
for each real number, $x$. From the inclusion
$$
\{n\in \N; v(n) < x\} \subset \{\alpha \in \Omega; \tilde{v}(\alpha)
\le x\}
$$
we obtain
$$
\text{cl}\;(\{n\in \N; v(n) < x\}) \subset \{\alpha \in \Omega;
\tilde{v}(\alpha) \le x\}.
$$
Furthermore (\ref{e2}) yields
$$
\mu^\ast(\{n\in \N; v(n) < x\}) \le F(x)
$$
for every real number $x$. On the other hand
$$
\N \setminus \{n\in \N; v(n) < x\} = \{n\in \N; v(n) \ge x\},
$$
therefore
$$
\text{cl}\;(\N \setminus \{n\in \N; v(n) < x\}) \subset\{\alpha \in \Omega;
\tilde{v}(\alpha) \ge x\}.
$$
Hence
$$
\mu^\ast(\N \setminus \{n\in \N; v(n) < x\}) \le 1-F(x),
$$
and so the set $\{n\in \N; v(n) < x\}$ is Buck measurable and its
measure density is $F(x)$.

The implication (ii) $\Rightarrow$ (iii) is trivial.

(iii) $\Rightarrow$ (i).

 Clearly
$$
\{n\in \N; v(n) < x\} \subset \{\alpha \in \Omega; \tilde{v}(\alpha)
\le x\})
,$$
and so $F(x) \le P(\tilde{v} \le x)$. On the other hand
$$
\{\alpha \in \Omega; \tilde{v}(\alpha) < x\}) \subset \text{cl}(\{n\in \N;
v(n) \le x\} )
$$
for $\varepsilon >0$. This yields $F(x) \le P(\tilde{v} \le x) \le
F(x+\varepsilon)$ for $\varepsilon>0$. For $\varepsilon \to 0^+$ we
obtain the assertion from the continuity of $F$.\qed

\section{Independence and measurability}

From our definitions in previous sections we immediately derive:
\begin{proposition}
 The Buck measurable sequences
$\{v_1(n)\},\{v_2(n)\}, \dots, \{v_r(n)\} $ are independent if and
only if for every $x_1,\dots,x_r \in \mathbb{R}$  we have
$$
\mu\Big(\bigcap_{j=1}^r \{n\in \N; v_j(n) < x_j\}\Big) =
\prod_{j=1}^r\mu(\{n\in \N; v_j(n) < x_j\}).
$$
\end{proposition}

\begin{example} We come back to Example \ref{vdc}. Consider
the sequences $\{Q_k^{(j)}\}$ given such that $Q_0^{(j)}=1,
j=1,\dots ,r$ and $Q_k^{(j)}/
Q_{k+1}^{(j)}$ for $j=1,\dots ,r$ and
$k=0,1,2 \dots$. Let $Q_k^{(j)}, Q_k^{(j_1)}$ be relatively prime
for $j \neq j_1$. Denote by $\{\gamma_j(n)\}$ the van der Corput
sequence with base $Q_k^{(j)}$ for $j=1,\dots ,r$. Then these
sequences are independent (see [IPT]).
\end{example}

\begin{theorem}
\label{indep} { Let $\{v_1(n)\},\{v_2(n)\}, \dots ,\{v_k(n)\}$ be
independent Buck measurable $p$-continuous sequences with continuous
Buck distribution functions $F_j, j=1\dots,k$. Then the random
variables $\tilde{v_1}, \dots, \tilde{v_k}$ are independent.}
\end{theorem}

{\bf Proof.} For $x_1,\dots x_k \in \mathbb{R}$  we have
$$
\{\alpha \in \Omega; \tilde{v_1}(\alpha) < x_1 ,\dots ,
\tilde{v_k}(\alpha) < x_k\} \subset
$$
$$
\subset \text{cl}(\{n\in \N; v_1(n) \le x_1 , \dots , v_k(n) \le
x_k\}).
$$
Thus $P(\tilde{v_1}<x_1 , \dots , \tilde{v_k}<x_k) \le
F_1(x) \dots F_k(x_k),$ and so from the above theorem we get
$P(\tilde{v_1}<x ,  \dots , \tilde{v_k}<_k)\le
P(\tilde{v_1}<x) \dots P(\tilde{v_k}<x_k)$.

On the other hand we have
$$
 P(\tilde{v_1}\le x_1) \dots P(\tilde{v_k}\le x_k)=
$$
$$
= \mu(\{n\in \N; v_1(n) \le x\}) \dots \mu(\{n\in \N; v_k(n) \le
 x_k\})=
$$
$$
=P(\text{cl}(\{n\in \N; v_1(n) \le x_1 , \dots , v_k(n) \le x_k\})
\le P(\tilde{v_1}\le x_1 ,  \dots , \tilde{v_k}\le x_k).
$$
 \qed

Let $F_1,...,F_k$ be non-decreasing functions defined on $\mathbb{R}$, ($k$ is a fixed positive integer). A set $B \subset \R^k$ is
called {\it Jordan Stieltjes measurable} with respect to the functions
$F_1,...,F_k$ if and only the Riemann Stieltjes integral
$$
\int \int ...\int \mathcal{X}_B dF_1...dF_k
$$
exists; $\mathcal{X}_B$ denoting the indicator function of $B$.

\begin{theorem}
\label{jordan} Let $\{v_1(n)\}, ..., \{v_k(n)\}$ be independent Buck
measurable $p$-continuous sequences with continuous Buck
distribution functions $F_1,...,F_k$. Suppose that a set $B \subset
\R^k$ {\it is Jordan Stieltjes measurable} with respect to the functions
$F_1,...,F_k$. Then the set $\{n \in \N; (v_1(n), ...,v_k(n)) \in
B\}$ is Buck measurable and its Buck measure density is
\begin{equation}
\label{integral} \int \int ...\int \mathcal{X}_B dF_1...dF_k.
\end{equation}
\end{theorem}

{\bf Proof.} If $B= [a_1,b_1] \times ...\times [a_k,b_k]$ is a
cylinder set then (\ref{integral}) follows directly from
independence of $\{v_1(n)\}, ..., \{v_k(n)\}$ and Theorem
\ref{thm1}. Proposition \ref{propI} then implies the assertion. \qed

\section{Integral and mean value}

Let $h: \Omega \to (-\infty, \infty)$ be a continuous function.
Since $\Omega$ is a compact space, it is uniformly continuous.
Consider  $m \in \N$. To the function $h$ we can associate a
periodic function $h_m$ with period $m$ in following way:
$$
 \alpha \in s+m\Omega  \Longleftrightarrow  h_m(\alpha) = h(s).
$$
Clearly,
\begin{equation}
\label{intm} \int h_m dP = \frac{1}{m} \sum_{s=0}^{m-1} h(s).
\end{equation}
Clearly $\lim_{N \to \infty}\mathfrak{d}(N!,0)=0,$ and so
uniform continuity of $h$ implies that $h_{N!}$ converges uniformly
to $h$. From (\ref{intm}) we obtain
\begin{equation}
\label{intinf} \int h dP = \lim_{N \to \infty} \frac{1}{N!}
\sum_{s=0}^{N!-1} h(s).
\end{equation}
The function $h$ restricted on $\N$ is $p$-continuous. Thus there
exists the limit $\lim_{m\to
\infty}\frac{1}{m}\sum_{s=0}^{m-1} h(s)$. From (\ref{intinf}) we
conclude
\begin{equation}
\label{INT} \int h dP = \lim_{m \to
\infty}\frac{1}{m}\sum_{s=0}^{m-1}h(s).
\end{equation}

\begin{remark}
If the random variable $\tilde{v}$ has a continuous distribution
function $F$ then
$$
\int \tilde{v} dP = \int_{-\infty}^{\infty} x dF(x) = E(v).$$
\end{remark}

The central limit theorem immediately yields:
\begin{proposition} { Let $\{v_k(n)\}, k=1,2,3 \dots$ be a sequence of
$p$-continuous sequences such that for every $k=1,2,3, \dots$ the
sequences $\{v_j(n)\}, j=1\dots k$ are independent and have the same
continuous Buck distribution function. Then for every $x \in
\mathbb{R}$ we have
$$
\lim_{k \to \infty}\mu\Big(\Big\{ n \in \N; \frac{v_1(n)+\dots
+v_k(n) - kE}{\sqrt{k}D} \le x \Big\}\Big) =
\frac{1}{\sqrt{2\pi}}\int_{-\infty}^x e^{\frac{-t^2}{2}} dt.
$$}
\end{proposition}

We conclude this section with the following metric result:

\begin{theorem} Let $v_k, k=1,2,3,\dots$ be a system of  independent $p$-continuous uniformly
distributed sequences. Then the sequence $\{\tilde{v}_n(\alpha)\}$
is uniformly distributed for almost all $\alpha \in \Omega$.
\end{theorem}

{\bf Proof.} Denote
$$
S_N(h, \alpha)=\frac{1}{N}\sum_{n=1}^N e^{2\pi i h
\tilde{v}_n(\alpha)}
$$
for $h \in \mathbb{Z}\setminus\{0\}$ and $\alpha \in
\Omega$. Put
$$
A_h =\{\alpha \in \Omega; \lim_{N \to \infty} S_N(h, \alpha)=0\}
$$
for $h\neq 0$. For every $n \in \N$ we have $E(e^{2\pi i h
\tilde{v}_n})=0$.  Therefore the strong low of large numbers implies
that $P(A_h)=1$. Thus $P(\cap_{h \neq 0} A_h)=1$ and the assertion
follows. \qed

\section{Weak Buck measurability}
\begin{definition}
\label{WPC} Let $v =\{v(n)\}$ be a real valued sequence. We say that
$v$ is {\it weakly polyadicly continuous} if and only if for every
$\varepsilon>0, \delta>0$ there exists a set $A \in \D_\mu$ with
$\mu(A)<\delta$ such that
$$
n_1 \equiv n_2 \pmod m \Rightarrow |v(n_1)-v(n_2)| < \varepsilon
$$
for all $n_1, n_2 \in \N \setminus A$ and a suitable $m\in \mathbb{N}$.
\end{definition}

Our aim is to prove the following equivalence :

\begin{theorem}
\label{EQ} A bounded sequence of real numbers is weakly Buck
measurable if and only if it is weakly polyadicly continuous.
\end{theorem}

We start by the proof of the first implication. We recall the following
notion:

\begin{definition} A real valued sequence $v$ is called {\it almost
polyadicly continuous} if and only if for each $\delta>0$ there
exists a set $A \in \D_\mu$ with $\mu(A)<\delta$ such that $v$ is
polyadicly continuous on the set $\N \setminus A$.
\end{definition}

Directly from the definition we get

\begin{proposition}
\label{indicator} A set $S\subset \N$ is  Buck measurable if and
only if its indicator function $\mathcal{X}_S$ is almost polyadicly
continuous.
\end{proposition}

\begin{proposition}
\label{linkom} If $v_1,v_2$ are two almost polyadicly continuous
sequences and $c_1, c_2$ are real numbers then the sequence
$c_1v_1+c_2v_2$ is almost polyadicly continuous.
\end{proposition}

\begin{proposition}
\label{APR} If $v$ is a real valued sequence such that for each
$\varepsilon>0$ there exists an almost polyadicly sequence $v_0$ that
$|v(n)-v_0(n)|< \varepsilon$ for $n \in \N$ then $v$ is weakly
polyadicly continuous.
\end{proposition}

\begin{proposition} Each bounded weakly Buck measurable sequence is weakly
polyadicly continuous.
\end{proposition}
{\bf Proof.} Let $v$ be a weakly Buck measurable sequence of elements
in the interval $[a,b], a<b$. Consider $\varepsilon>0$. Then there
exists a partition $x_0, \dots, x_k$  of $[a,b]$ such that the
sets
$$
S_i =\{n \in \N; v(n)
\in [x_i, x_{i+1})\}, i=0,\dots, k-2
$$
and
$$
S_{k-1}= \{n \in \N; v(n) \in [x_{k-1}, b]\}
$$
are Buck measurable and $x_{i+1}-x_i < \varepsilon$. Then the
sequence
$$
v_0(n)=\sum_{i=0}^{k-1} x_i\mathcal{X}_{S_i}(n), n \in \N
$$
is almost polyadicly continuous and $|v_0(n)-v(n)| < \varepsilon$.
The assertion follows from Proposition \ref{APR}.\\

Now we prove the second implication.

If $v = \{v(n)\}$ is a real valued sequence and $k=\{k_n\}$ is a
sequence of positive integers then we shall denote
$v(k)=\{v(k_n)\}$.

\begin{proposition}
\label{meskrit} A set $S \subset \N$ is Buck measurable if and only
if for each sequence of positive integers $k$ the sequence
$\mathcal{X}_S(k)$ has a mean value and in this case
$$
\mu(S)=E(\mathcal{X}_S(k)).
$$
\end{proposition}
This proposition is an easy reformulation of Theorem 7 in [PAS3]
page 51 or Theorem 50 in [PAS5] page 113.
\begin{proposition}
\label{pp34} If $v$  is a bounded weakly polyadicly continuous
sequence then it has a mean value and for each sequence of positive
integers $k$ which is uniformly distributed in $\Z$ we have
$$
E(v(k))=E(v).
$$
\end{proposition}

{\bf Proof.} Consider $\delta>0, \varepsilon >0$. Let $A, m$ be as
in Definition \ref{WPC}. Suppose that $r_1, \dots, r_s$ is the
maximal finite sequence of elements of $\N \setminus A$ incongruent
modulo $m$ and $r_{s+1}, \dots, r_m$ its completion with respect to a complete
residue system modulo $m$. Define the periodic sequence $v_m(n) =
v(r_j)$ if and only if $n \equiv r_j \pmod{m}$ for $j=1, \dots, m$
and $n \in \N$. Then for each $n \in \N \setminus A$ we have
\begin{equation} \label{approx}
|v_m(n) - v(n)| < \varepsilon.
\end{equation}
For $N =1,2,3, \dots$ we obtain
$$
E_N(v_m(k))-E_N(v(k))= \frac{1}{N}\sum_{n=1}^N (v_m(k_n)-v(k_n))=
$$
$$
\frac{1}{N}\sum_{n\le N, k_n \in A}
(v_m(k_n)-v(k_n))+\frac{1}{N}\sum_{n\le N, k_n \not\in A}
(v_m(k_n)-v(k_n)).
$$
And so from (\ref{approx}) we device
$$
|E_N(v_m(k))-E_N(v(k))| < 2HE_N(\mathcal{X}_A(k))+\varepsilon
$$
where $H$ is upper bound of $\{|v(n)|\}$. If $k$ is uniformly
distributed in $\Z$ we get for $N \to \infty$
\begin{equation}
\label{wpc1} |E(v_m)-\overline{E}(v(k))| < 2H\delta+\varepsilon
\end{equation}
and
\begin{equation}
\label{wpc2} |E(v_m)-\underline{E}(v(k))| < 2H\delta+\varepsilon.
\end{equation}
Therefore
$$
\overline{E}(v(k)) - \underline{E}(v(k)) <4H\delta+2\varepsilon.
$$
Since $\delta, \varepsilon$ are arbitrary we have
$\overline{E}(v(k))=\underline{E}(v(k))=E(v(k))$. If in the
inequalities (\ref{wpc1}) and (\ref{wpc2}) we substitute  the
sequence $\{n\}$ instead of $k$ we conclude $E(v)=E(v(k))$. \qed

\begin{proposition}
\label{transfor} If $v$ is a weakly polyadicly continuous sequence of
elements in  $[a,b]$ and $f$ is a continuous real function
defined on this interval then the sequence $f(v)$ is weakly
polyadicly continuous, too.
\end{proposition}

{\bf Proof.} The assertion follows immediately from the fact that a
continuous function on a compact interval is uniformly continuous. \qed

\begin{proposition} Each bounded weakly polyadic continuous real
valued sequence is weakly Buck measurable.
\end{proposition}

{\bf Proof.} Let $v$ be a weakly polyadic continuous real valued
sequence of elements in $[a,b].$ Then for
every continuous function $f$ defined on $[a,b]$ the sequence $f(v)$
has a mean value and for every sequence of positive integers $k$
uniformly distributed in $\Z$ we have
$$
E(f(v))=E(f(v(k)).
$$
We define a positive linear functional
$$
\Phi(f)=E(f(v))
$$
on the linear space of all continuous real functions defined on
$[a,b]$ such that $\Phi(1)=1$.

Thus Riesz representation theorem provides that a non decreasing
function $F$ exists such $F(a)=0, F(b)=1$ and
\begin{equation}
\label{riesz} E(f(v(k))=\Phi(f)= \int_a^b f(x)dF(x)
\end{equation}
holds for each sequence $k$ uniformly distributed in $\Z$. If the
function $F$ is continuous in $x_0$ then by Proposition \ref{APP} we
can construct for every $\varepsilon > 0$ two continuous functions
$f_1,f_2$ defined on $[a,b]$ satistying
$$
\int_a^b(f_2(x) - f_1(x)) dF(x) < \varepsilon
$$
and
$$
f_1 \le \mathcal{X}_{[0,x_0)} \le f_2.
$$
Hence for each sequence $k$ uniformly distributed in $\Z$ we have
$$
E(\mathcal{X}_{[0,x_0)}(v(k)))=F(x_0).
$$
Proposition \ref{meskrit} implies that  the set $\{n \in \N;
v(n)<x_0\}$ is Buck measurable. Since every non-decreasing function has at
most a countable set of discontinuities, the proof is complete. \qed

Let $\mathcal{B}_\mu$ be the set of all bounded weakly measurable
sequences. Theorem \ref{weakmesur} implies

\begin{proposition}
\label{Ba} Define the norm
$$
||v|| = \sup\{|v(n)|; n \in \N\}
$$
for $v \in \mathcal{B}_\mu.$ Then $(\mathcal{B}_\mu, +, \cdot, ||
\cdot ||)$ is a Banach algebra.
\end{proposition}

\section{Statistical independence}

If $v=\{v(n)\}$ is a sequence and $g$ is a function defined on the
set containing the elements of $v$ then we denote by
$g(v)\;\text{the sequence}\;\{g(v(n))\}$. The following theorem relates $p$-continuous independent sequences
to the concept of uniform distribution in $\mathbb{Z}.$

\begin{theorem}
\label{sss} { Let $\{v_1(n)\},\dots,\{v_k(n)\}$ be $p$-continuous
independent sequences. Then for arbitrary functions $g_1,\dots,g_k$
continuous on the real line we have
$$
 E\Big(\prod_{j=1}^k g_j(v_j(k)\Big)
= \prod_{j=1}^k E\Big(g_j(v_j(k))\Big)
$$
for each sequence $k = \{k_n\}$ which is uniformly distributed in
$\Z$.}
\end{theorem}

{\bf Proof.} If $\{v(n)\}$ is a $p$-continuous function, then
it is bounded. Every continuous function $g$ defined on the real line is
uniformly continuous on the interval $[b_1, b_2]$ where $b_1$ is a
lower bound of  the sequence $\{v(n)\}$ and $b_2$ its upper bound.
Thus the sequence  $\{g(v(n))\}$ is $p$-continuous, too.

 Let us consider $\{v_1(n)\},\dots,\{v_k(n)\}$ -
polyadicly continuous independent sequences. Then the random
variables $\tilde{v_1},\dots,\tilde{v_k}$ are independent and
so the random variables $g_1(\tilde{v_1}),\dots,g_k(\tilde{v_k})$
are independent, too. Thus
$$
E(g_1(v_1)\dots g_k(v_k)) = E(g_1(v_1))\dots E(g_k(v_k))
$$
and the assertion follows from Proposition \ref{pp34}.\\

Theorem \ref{GENIND} and Proposition \ref{pp34} imply
\begin{theorem}
{ Let $\{v_1(n)\},\dots,\{v_k(n)\}$ Buck measurable independent
sequences having continuous Buck distribution functions. Then for
arbitrary functions $g_1,\dots,g_k$ continuous on the real line
$$
 E\Big(\prod_{j=1}^k g_j(v_j(k)\Big)
= \prod_{j=1}^k E\Big(g_j(v_j(k))\Big)
$$
for each sequence $k = \{k_n\}$ which is uniformly distributed in
$\Z$.}
\end{theorem}

Adresses of authors :

Milan Pa\v st\'eka,  Department of Mathematics and Informatics,
Faculty of education, Univesity of Trnava, Priemyseln\'a 4, P. O.
BOX 9,  918 43  Trnava, Slovakia

Robert Tichy, Institut fur Analysis und Zahlentheorie, Technische
Universitat Graz, Steyrergasse 30 A-8010, Graz, Austria


\begin{thebibliography}{TMF9}
 \addcontentsline{toc}{section}{\quad\ \ Bibliography}

\bibitem [BUC]{BUC}
{{\sc Buck, R., C.,}} \textsl{The measure theoretic approach to
density},
 Amer. J. Math \textbf{68}, 1946, 560--580

\bibitem [CQ]{CQ}
{{\sc Choimet, D., Queffelet, H.,}} \textsc{Analyse mathematique,
grandes theoremes du vingtieme siecle} Calvage Mounet, Paris, (2009)


\bibitem [D-T]{D-T}
{{\sc Drmota, M., Tichy, R. F.,}} \textsl{ Sequences, Discrepancies
and Applications, Springer, Berlin Heidelberg},
 Springer, Berlin Heidelberg,
1997


\bibitem [Er]{Er}
{{\sc Erd\H{o}s, P.,}} \textsl{ On the density of some sequences of
numbers  II,} Journal of the London  Math. Soc. 12, (1937), 7 - 11

\bibitem [G]{G}
{{\sc Grekos, G.,}} \textsl{On various definitions of density (a
survey)}, Tatra Mt. Math. Publ.,  31, 2005, 17--27
\bibitem [G1]{G1}
{{\sc Grekos, G.,}} \textsl{The density set (a survey)}, Tatra Mt.
Math. Publ.,  31, 2005, 103--111
\bibitem [GST] {GST}
{{\sc Grabner, P. J., Strauch, O., Tichy, R. F.,}}
\textsl{Lp-discrepancy and statistical independence of sequences, }
Czechoslovak Mathematical Journal, Vol.49(1999), No.1,97�110

\bibitem [GT] {GT}
{{\sc Grabner, P. J., Tichy, R. F.,}} \textsl{Remark on statistical
independence of sequences, } Mathematica Slovaca,
Vol.44 (1994), No.1,91--94


\bibitem [IPT]{IPT}
{{\sc  Iaco, M. R., Pasteka, M.,  Tichy R., F.,}} \textsl{Measure
density for set decompositions and uniform distribution} Rend. Circ.
Math. Palermo (2) , 64, No. 2,  2015 , 323 -- 339


\bibitem [K-N]{K-N}
{{\sc  Kuipers, L.,  Niederreiter, H., }} \textsl{ Uniform
distribution of Sequences}, John Wiley and Sons, N.Y. London, Sydney
Toronto, 1974


\bibitem [N]{N}
{{\sc Novoselov, E. V.,}} \textsl{Topological theory of polyadic
numbers}, Trudy Tbilis. Mat. Inst. 27, 1960, 61 -- 69, (in russian)
\bibitem [N1]{N1}
{{\sc Novoselov, E. V.,}} \textsl{New method in probabilistic
number theory}, Doklady akademii nauk. ser. matem. No. 2, 28, 1964,
307 -- 364,  in russian

\bibitem [NAR]{NAR}
 {{\sc Narkiewicz, W. , }}
 \textsl{Teoria liczb}, (in polish)
 PWN, Warszawa, 1991

\bibitem [NIV]{NIV}
 {{\sc Niven, I.,}}
 \textsl{Uniform distribution of sequences of integers},
 Trans. Amer. Math. Soc. 98, 52 -- 61


\bibitem [PAS]{PAS}
{{\sc  Pa\v st\'eka, M.,}} \textsl{Some properties of Buck's measure
density\index{measure density}},
  Math. Slovaca 42, no. 1,  1992, 15--32

\bibitem [PAS2]{PAS2}
{{\sc  Pa\v st\'eka, M.,}} \textsl{ Remarks on one type of uniform
distribution} Unif. Distrib. Theory 2, No. 1, 2007, 79--92

\bibitem [PAS3]{PAS3}
{{\sc  Pa\v st\'eka, M.,}} \textsl{On four approaches to density}
Spectrum Slovakia 3. Frankfurt am Main: Peter Lang; Bratislava:
VEDA, Publishing House of the Slovak Academy of Sciences, 2014


\bibitem [PAS4]{PAS4}
{{\sc  Pa\v st\'eka, M.,}} \textsl{Remarks on Buck's measure
density}. Tatra Mt. Math. Publ. 3, 1993, 191--200


\bibitem [PAS5]{PAS5}
{{\sc  Pa\v st\'eka, M.,}} \textsl{Density and related topics},
Veda, Bratislava, Academia, Praha, 2017.
\bibitem [P-T]{P-T}
{{\sc Pa\v st\'eka, M., Tichy, R.,}} \textsl{A note on the correlation
coefficient of arithmetic functions }, Acta Acad. Paed. Agriensis,
Sectio Mathematicae  30, 2003, 109--114

\bibitem [Ra] {Ra}
{{ \sc Rauzy, G.,}} \textsl{Proprietes statistique de suites
arithmetiques.} Le Mathematicien. No. li. Collection SUP, Presses
Universitaires de France, Paris, 1976

\bibitem [Sch]{Sch}
{{\sc Schoenberg, I.,}} \textsl{\"Uber die asymptotische
Verteilung reeller Zahlen mod 1.},
  Math. Z., 1928, 171 -- 199


\bibitem [SP]{SP}
{{\sc Strauch, O. Porubsk\'y, \v S.,}}\textsl{  Distribution of
Sequences a Sampler, Peter Lang, SAV, Frankfurt am Main}, Peter
Lang, SAV, Frankfurt am Main, 2005

\bibitem [WEY]{WEY}
 {{\sc Weyl, H.}},
 \textsl{\"Uber die Gleichverteilung von Zahlen mod. Eins},
 Math. Ann, 77, 1916, 313--352

\end{thebibliography}
\end{document}